\documentclass[11pt]{amsart}
\usepackage{amscd}
\usepackage{curves}
\usepackage{epsfig}
\usepackage[pass]{geometry}
\usepackage{mathrsfs}
\usepackage{pstricks}
\usepackage{pstricks-add}
\usepackage{pst-grad}
\usepackage{pst-plot}
\usepackage{auto-pst-pdf}
\usepackage{verbatim}
\usepackage{latexsym,eucal,amsfonts,amssymb,amsmath,graphicx}
\numberwithin{equation}{section}
\newtheorem{theorem}{Theorem}[section]
\newtheorem{lemma}[theorem]{Lemma}
\newtheorem{prop}[theorem]{Proposition}
\newtheorem{cor}[theorem]{Corollary}

\newtheorem{remark}[theorem]{Remark}
\def \bpf {\begin{proof}}
\def \epf {\end{proof}}
\def \beq {\begin{equation*}}
\def \eeq {\end{equation*}}

\def \mcd {{\mathscr D}}

\def \mcs {{\mathscr S}}

\def \mcu {{\mathscr U}}

\def \mcx {{\mathscr X}}
\def \mbr {{\mathbb R}}

\def \id {\operatorname{Id}}
\def \comp {\operatorname{comp}}
\def \loc {\operatorname{loc}}

\def \diag{\textrm{Diag}}

\def \La {\Lambda}    
   
\def \lap {\Delta}
\def \p {\partial}
\def \eps {\epsilon}
\def \det {\text{det}}

\def \ha {\frac{1}{2}}

\def \beqq {\begin{equation}}
\def \eeqq {\end{equation}}

\def \WF {\text{WF}}

\begin{document}
\title[Nonlinear responses from the interactions of two waves]{Nonlinear responses from the interaction of two progressing waves at an interface}
\author{Maarten de Hoop}
\address{Maarten de Hoop
\newline
\indent Computational and Applied Mathematics and Earth Science, Rice University}
\email{mdehoop@rice.edu}

\author{Gunther Uhlmann} 
\address{Gunther Uhlmann
\newline
\indent Department of Mathematics, University of Washington,  
\newline
\indent Institute for Advanced Study, the Hong Kong University of Science and Technology,
\newline
\indent and Department of Mathematics, University of Helsinki}
\email{gunther@math.washington.edu}

\author{Yiran Wang}
\address{Yiran Wang
\newline
\indent Department of Mathematics, University of Washington 
}
\email{wangy257@math.washington.edu}
\begin{abstract} 
For scalar semilinear wave equations, we analyze the interaction of two (distorted) plane waves at an interface between media of different nonlinear properties. We show that new waves are generated from the nonlinear interactions, which might be responsible for the observed nonlinear effects in applications. Also, we show that the incident waves and the nonlinear responses determine the location of the interface and some information of the nonlinear properties of the media. In particular, for the case of a jump discontinuity at the interface, we can determine the magnitude of the jump.
\end{abstract}

\maketitle

\section{Introduction}
Let $g$ be a smooth Riemannian metric on $\mbr^3$. In local coordinates $x = (x^1, x^2, x^3)$, the (positive) Laplace-Beltrami operator is given by
\beq
\lap_g = - \frac{1}{\sqrt{\det g}} \sum_{i, j = 1}^3\frac{\p}{\p x^i}(\sqrt{\det g} g^{ij} \frac{\p}{\p x^j}).
\eeq
We shall work with the associated wave operator 
\beq
P = \p_t^2 + \lap_g.
\eeq
However, one can consider $P$ with lower order perturbations to which the results of this work apply as well. For example, one can consider wave operators with variable sound speed and density
\beq
\widetilde P = \p_t^2 - c^2(t, x) \rho(t, x) \nabla \cdot(\frac{1}{\rho(t, x)}\nabla u),
\eeq
where $c(t, x)$ is the sound speed and $\rho(t, x)$ is the density of the medium. 

Consider the following semilinear wave equation 
\beqq\label{eqsem0}
\begin{gathered}
P u(t, x)  + a(t, x) u^2(t, x) = 0, \ \ \text{ in }  (-\infty \times T) \times \mbr^3,\\
u(t, x) = u_0(t, x), \ \ \text{ in } (-\infty, 0)\times \mbr^3,
\end{gathered}
\eeqq
with $T> 0.$ Suppose that the incident wave $u_0$ consists of progressing plane waves with conormal singularities to two characteristic surfaces $S_1$ and $S_2$ for $P$ 
which do not intersect for $t<0$. When $a$ is smooth and the spatial dimension is two, the interaction of waves was studied in Bony \cite{Bo}, Melrose-Ritter \cite{MR1} and others. In particular, as a special case of \cite[Theorem 1]{MR1}, we know that the solution is conormal to $S_1$ and $S_2$ after the interaction and no new wave is produced. Melrose and Ritter \cite[Theorem 2]{MR1} showed that the interaction of three progressing waves could generate new waves. Explicit examples when the new waves are indeed produced have been constructed by various authors; see Rauch-Reed \cite{RR} and the text book by Beals \cite{Bea}. For $a$ smooth and  spatial dimension three, such phenomena have also been analyzed and the newly generated waves have played an important role in the inverse problem for nonlinear hyperbolic equations in \cite{KLU, KLU1, LUW}.

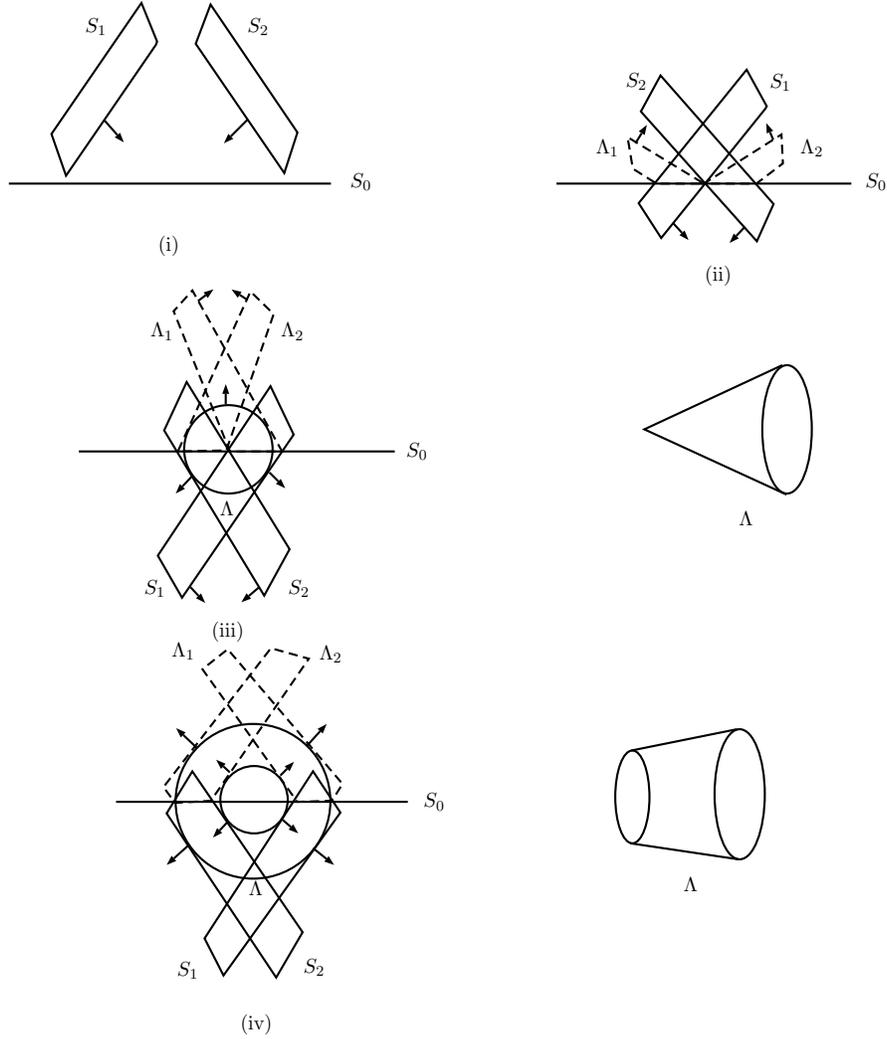
\begin{figure}[htbp]
\centering
\scalebox{0.7} 
{
\begin{pspicture}(0,-2.720867)(16.596666,2.720867)
\psline[linecolor=black, linewidth=0.04](0.0,-0.7542002)(6.116667,-0.7542002)
\psline[linecolor=black, linewidth=0.04](10.4,-0.7542002)(16.0,-0.7542002)
\pspolygon[linecolor=black, linewidth=0.04](14.011491,1.4107454)(11.960923,-1.1262925)(12.38875,-1.788689)(14.367307,0.66958755)(14.393362,0.7160833)
\pspolygon[linecolor=black, linewidth=0.04](12.375893,1.3010976)(12.002668,0.6157823)(14.206059,-1.8506409)(14.523111,-1.1423936)
\pspolygon[linecolor=black, linewidth=0.04](2.5173829,2.678812)(0.8087431,0.18918175)(1.0850531,-0.6038922)(2.7936928,1.885738)(2.8145497,1.9462168)
\pspolygon[linecolor=black, linewidth=0.04](3.8304713,2.6488242)(3.5485063,1.9036024)(5.23041,-0.55371076)(5.4801726,0.21764219)
\pspolygon[linecolor=black, linewidth=0.04, linestyle=dashed, dash=0.17638889cm 0.10583334cm](12.3,-0.7542002)(11.85,-0.4708669)(11.75,0.12913309)(13.2,-0.7542002)
\pspolygon[linecolor=black, linewidth=0.04, linestyle=dashed, dash=0.17638889cm 0.10583334cm](14.1973505,-0.7542002)(14.7,-0.3708669)(14.6725445,0.1791331)(13.2,-0.7542002)
\rput[bl](6.483333,-0.88753355){$S_0$}
\rput[bl](1.45,2.0624664){$S_1$}
\rput[bl](4.5333333,2.0624664){$S_2$}
\rput[bl](11.716666,1.0291331){$S_2$}
\rput[bl](14.45,1.0124664){$S_1$}
\rput[bl](16.266666,-0.88753355){$S_0$}
\rput[bl](2.8333333,-2.120867){(i)}
\rput[bl](13.216666,-2.7042003){(ii)}
\psline[linecolor=black, linewidth=0.04, arrowsize=0.05291667cm 2.0,arrowlength=1.4,arrowinset=0.0]{->}(1.8166667,0.44579977)(2.1833334,0.029133098)
\psline[linecolor=black, linewidth=0.04, arrowsize=0.05291667cm 2.0,arrowlength=1.4,arrowinset=0.0]{->}(4.5333333,0.46246642)(4.0833335,0.012466431)
\psline[linecolor=black, linewidth=0.04, arrowsize=0.05291667cm 2.0,arrowlength=1.4,arrowinset=0.0]{->}(12.607654,-1.4970727)(12.909013,-1.8279943)
\psline[linecolor=black, linewidth=0.04, arrowsize=0.05291667cm 2.0,arrowlength=1.4,arrowinset=0.0]{->}(13.980293,-1.5748265)(13.686373,-1.8835739)
\psline[linecolor=black, linewidth=0.04, arrowsize=0.05291667cm 2.0,arrowlength=1.4,arrowinset=0.0]{->}(11.9,0.012466431)(12.116667,0.36246642)
\psline[linecolor=black, linewidth=0.04, arrowsize=0.05291667cm 2.0,arrowlength=1.4,arrowinset=0.0]{->}(14.516666,0.0791331)(14.383333,0.39579976)
\rput[bl](11.15,-0.23753357){$\Lambda_1$}
\rput[bl](15.033334,-0.23753357){$\Lambda_2$}
\end{pspicture}
}
\scalebox{0.7} 
{
\begin{pspicture}(0,-3.3467042)(13.94,3.3467042)
\psline[linecolor=black, linewidth=0.04](0.0,0.27329576)(6.0,0.27329576)
\pspolygon[linecolor=black, linewidth=0.04](3.64,1.5132958)(1.5,-1.7067043)(1.96,-2.5067043)(4.08,0.59329575)
\pspolygon[linecolor=black, linewidth=0.04](1.62,0.67329574)(3.52,-2.4667044)(4.00254,-1.58531)(2.0477226,1.592554)
\pscircle[linecolor=black, linewidth=0.04, dimen=outer](2.84,0.31329575){0.86}
\pspolygon[linecolor=black, linewidth=0.04, linestyle=dashed, dash=0.17638889cm 0.10583334cm](1.88,0.27329576)(3.26,3.3132958)(3.7,2.8834162)(2.8,0.2929343)
\pspolygon[linecolor=black, linewidth=0.04, linestyle=dashed, dash=0.17638889cm 0.10583334cm](3.86,0.27329576)(2.1565385,3.3332958)(1.8,2.938457)(2.8698077,0.2722635)
\rput[bl](1.24,-2.4867043){$S_1$}
\rput[bl](3.98,-2.4867043){$S_2$}
\rput[bl](6.22,0.11329575){$S_0$}
\rput[bl](2.68,-0.94670427){$\La$}
\psellipse[linecolor=black, linewidth=0.04, dimen=outer](13.45,0.6932957)(0.49,1.24)
\psline[linecolor=black, linewidth=0.04](13.36,1.9132957)(10.74,0.6932957)(13.44,-0.54670423)
\rput[bl](12.54,-1.1267042){$\La$}
\rput[bl](2.52,-3.3467042){(iii)}
\psline[linecolor=black, linewidth=0.04, arrowsize=0.05291667cm 2.0,arrowlength=1.4,arrowinset=0.0]{->}(2.12,-2.3067043)(2.4,-2.6067042)
\psline[linecolor=black, linewidth=0.04, arrowsize=0.05291667cm 2.0,arrowlength=1.4,arrowinset=0.0]{->}(3.4334826,-2.29966)(3.0865173,-2.5937486)
\psline[linecolor=black, linewidth=0.04, arrowsize=0.05291667cm 2.0,arrowlength=1.4,arrowinset=0.0]{->}(2.262568,3.1319985)(2.557432,3.3545928)
\psline[linecolor=black, linewidth=0.04, arrowsize=0.05291667cm 2.0,arrowlength=1.4,arrowinset=0.0]{->}(3.2007983,3.15446)(2.8992016,3.3521316)
\psline[linecolor=black, linewidth=0.04, arrowsize=0.05291667cm 2.0,arrowlength=1.4,arrowinset=0.0]{->}(2.7894454,1.1730336)(2.7905548,1.5535579)
\psline[linecolor=black, linewidth=0.04, arrowsize=0.05291667cm 2.0,arrowlength=1.4,arrowinset=0.0]{->}(2.16,-0.18670425)(1.84,-0.52670425)
\psline[linecolor=black, linewidth=0.04, arrowsize=0.05291667cm 2.0,arrowlength=1.4,arrowinset=0.0]{->}(3.62,-0.12670426)(3.92,-0.42670426)
\rput[bl](3.84,2.3932958){$\Lambda_2$}
\rput[bl](1.36,2.3932958){$\Lambda_1$}
\end{pspicture}
}
\scalebox{0.7} 
{
\begin{pspicture}(0,-3.705)(12.34,3.705)
\psline[linecolor=black, linewidth=0.04](0.0,0.715)(5.54,0.715)
\pspolygon[linecolor=black, linewidth=0.04](3.74,1.295)(1.68,-1.8814286)(2.04,-2.585)(4.26,0.5107143)
\pspolygon[linecolor=black, linewidth=0.04](0.96,0.495)(3.04,-2.625)(3.54254,-1.7636057)(1.4477224,1.2942582)
\pscircle[linecolor=black, linewidth=0.04, dimen=outer](2.6,0.725){1.49}
\rput[bl](1.16,-2.645){$S_1$}
\rput[bl](3.56,-2.585){$S_2$}
\rput[bl](5.84,0.555){$S_0$}
\rput[bl](2.52,-1.065){$\La$}
\psellipse[linecolor=black, linewidth=0.04, dimen=outer](11.843098,0.8587908)(0.49690142,1.2562091)
\rput[bl](10.77831,-0.985){$\La$}
\rput[bl](2.36,-3.705){(iv)}
\pscircle[linecolor=black, linewidth=0.04, dimen=outer](2.62,0.755){0.66}
\psellipse[linecolor=black, linewidth=0.04, dimen=outer](9.804789,0.8081373)(0.34478873,0.901634)
\psline[linecolor=black, linewidth=0.04](9.784507,1.6895097)(11.792395,2.0947385)
\psline[linecolor=black, linewidth=0.04](9.764225,-0.073235296)(11.812676,-0.37715685)
\pspolygon[linecolor=black, linewidth=0.04, linestyle=dashed, dash=0.17638889cm 0.10583334cm](1.12,0.695)(0.9,0.995)(2.94,3.635)(3.66,3.435)(1.8,0.735)
\pspolygon[linecolor=black, linewidth=0.04, linestyle=dashed, dash=0.17638889cm 0.10583334cm](4.1,0.735)(4.3,1.055)(2.12,3.615)(1.64,3.235)(3.44,0.715)
\psline[linecolor=black, linewidth=0.04, arrowsize=0.05291667cm 2.0,arrowlength=1.4,arrowinset=0.0]{->}(1.36,-0.125)(0.96,-0.485)
\psline[linecolor=black, linewidth=0.04, arrowsize=0.05291667cm 2.0,arrowlength=1.4,arrowinset=0.0]{->}(3.74,-0.165)(4.14,-0.485)
\psline[linecolor=black, linewidth=0.04, arrowsize=0.05291667cm 2.0,arrowlength=1.4,arrowinset=0.0]{->}(3.62,1.775)(4.0,2.195)
\psline[linecolor=black, linewidth=0.04, arrowsize=0.05291667cm 2.0,arrowlength=1.4,arrowinset=0.0]{->}(1.54,1.735)(1.14,2.115)
\psline[linecolor=black, linewidth=0.04, arrowsize=0.05291667cm 2.0,arrowlength=1.4,arrowinset=0.0]{->}(2.16,1.275)(1.9,1.515)
\psline[linecolor=black, linewidth=0.04, arrowsize=0.05291667cm 2.0,arrowlength=1.4,arrowinset=0.0]{->}(2.12,0.335)(1.84,0.035)
\psline[linecolor=black, linewidth=0.04, arrowsize=0.05291667cm 2.0,arrowlength=1.4,arrowinset=0.0]{->}(3.16,0.375)(3.46,0.115)
\psline[linecolor=black, linewidth=0.04, arrowsize=0.05291667cm 2.0,arrowlength=1.4,arrowinset=0.0]{->}(3.12,1.215)(3.38,1.475)
\rput[bl](1.06,3.415){$\Lambda_1$}
\rput[bl](3.86,3.355){$\Lambda_2$}
\end{pspicture}
}
\label{figinter3}
\caption{Evolution of two plane waves interacting at an interface. In
  Figure (i), $S_1, S_2$ represent the wave fronts (singular supports)
  of two progressing waves in $\mbr^3$ and $S_0$ represents the
  singular support of $a(t, x)$. The picture shows the projective view
  on a plane $\mbr^2$ before the wave meets. The arrows indicate the
  directions of the wave propagation. Figure (ii) shows the
  intersection of the two waves at $S_0$ before they meet
  together. The dashed surfaces represent the reflected waves. Figure
  (iii) illustrates various waves during the interaction of the two
  waves at $S_0$. The wave front of the newly generated wave is
  demonstrated by the disk denoted by $\La$. The figure to the right
  shows the wave front in $\mbr^3$ which is the surface of a
  cone. Figure (iv) shows the waves after the interaction is
  complete. The wave front $\La$ actually becomes the surface of a
  truncated cone in $\mbr^3$ (picture to the right).}
\end{figure}

In this work, we are interested in the interactions of two progressive waves at an interface of media with difference nonlinear properties. In particular, we assume that $a(t, x)$ has conormal singularities at  a co-dimension one submanifold $S_0$ (the interface) of $\mbr^4$ not characteristic for $P$. A useful example to keep in mind is $a(t, x) = a(x)$ conormal to some $Y\subset \mbr^3$ regarded as the interface. For example, $a(x)$ or its derivatives have jump discontinuities across $Y$. If $S_1, S_2$ and $S_0$ intersect in $t \in (0, T)$ for some $T>0$ small, we show in Theorem \ref{main} that a new wave is produced due to the nonlinear interactions; see Figure \ref{figinter3} for an illustration of this interaction. In some sense, the nonlinear coefficient $a(t, x)$ plays the role of the third wave in the result mentioned above. 

The main motivation of our analysis comes from the study of nonlinear
interaction of waves related to conormal discontinuities
(``interfaces'') in the nonlinearities of the elastic moduli in
sedimentary rocks. Nonlinear properties of such rocks are commonly
associated with material damage. Nonlinear properties of solids have
been extensively studied in the laboratory by Rollins, Taylor and
Todd \cite{RTT}, Johnson, Shankland, O'Connell and Albricht \cite{JSO}, Johnson and
Shankland \cite{JS}, and many others. In the context of this paper, we are
concerned with so-called fast nonlinear dynamics (Johnson and
McCall \cite{JM}). Traditionally, the nonlinear interaction, in the absence of
singularities in the nonlinearities of the elastic moduli, has been
studied using monochromatic incident waves aiming to observe the
generation of combined harmonics; for an early analysis, see Jones and
Kobett \cite{JK}. (The experimental counterpart to our problem in some sense is
the one of two incident non-collinear beams generating a new beam at
their difference frequency.) This is also the underlying principle in
the scalar-wave formulation -- which we consider here -- for
vibro-acoustography \cite{FG1, FG2} based on ultrasound-stimulated
acoustic emission. However, the use of transient incident waves and
the generation (emission) of a new transient wave that we analyze,
here, has so far not been considered in applications and
experiments~\footnote{P.A. Johnson, personal communication}. Indeed,
the generation of this wave opens new ways for nonlinear imaging in
Earth's subsurface, which we elucidate here in the form of an inverse
problem. Studying the interaction with conormal singularities in the
nonlinearities of the elastic moduli was motivated by the work of
Kuvshinov, Smit and Campman \cite{KSC}. In a forthcoming paper, we
extend the results of this paper pertaining to scalar waves to the
elastic case.

We consider in Section \ref{sec-inv} an inverse problem and we apply the results of the previous sections. 
We send two distorted plane waves concentrated along geodesics that meet at the interface. We observe the nonlinear response. We show that from this information we can determine the interface and the principal symbol of $a(t, x)$. In particular, in the case that $a(t, x)$ has a jump type singularity we can determine the magnitude of the jump. For a precise statement of the problem and the results see Theorem \ref{main1}. 

The paper is organized as follows. In Section \ref{sec-lin}, we review
the theory for linear wave equations and construct distorted plane
waves as in \cite{KLU}. We establish local well-posedness of the
nonlinear wave equation with a non-smooth nonlinear term in Section
\ref{secwell}. In Section \ref{sec-asymp}, we analyze the nonlinear
responses after the interactions. In Section \ref{sec-linres}, we compare the linear and nonlinear responses in case the linear operator $P$ also has conormal singularities. We demonstrate that the conic wave is a distinctive feature of the nonlinear response. Finally, in Section \ref{sec-inv} we formulate and study
the inverse problem. 

\section*{Acknowledgement} 
Maarten V. de Hoop acknowledges and sincerely thanks the Simons Foundation under the MATH$+$X program for financial support. He was                  
also partially supported by NSF under grant DMS-1559587 and the members of the Geo-Mathematical Imaging Group at Rice University. 
Gunther Uhlmann is partially supported by NSF, a Si-Yuan Professorship at HKUST and FiDiPro Professorship of the Academy of Finland.

\section{The linear wave equation and distorted plane waves}\label{sec-lin}
We know (e.g.\ from \cite{Bar}) that for the linear wave equation
\beq
Pv = (\p_t^2 + \lap_g) v = f,
\eeq
 there exists a fundamental solution (causal inverse) $Q$ such that $QP = \id$ on the space of distributions $\mcd'(\mbr^4)$. We review the structure of the Schwartz kernel of the causal inverse. 

In the following, we use $x = (x^i)_{i = 0}^3$ as the local coordinates of $\mbr^4$ with $x^0 = t$. The dual variables in the cotangent bundle are denoted by $\zeta = (\tau, \xi), \tau\in \mbr, \xi\in \mbr^3$. Let $p(x, \zeta) = -\tau^2 + |\xi|^2_{g^*}$ be the symbol of $P$, where $g^* = g^{-1} = (g^{ij})$ is the dual metric. We denote by $\Sigma_P = \{(x, \zeta)\in T^*\mbr^4: p(x, \zeta) = 0\}$  the characteristic set for $P$ and $\Sigma_{P, x} \doteq \Sigma_P\cap T_x^*\mbr^4$ for any $x\in \mbr^4.$ The Hamilton vector field of $p(x, \zeta)$ is denoted by $H_p$ and in local coordinates
\beq
H_p = \sum_{i = 0}^3( \frac{\p p}{\p \zeta_i}\frac{\p }{\p x^i} - \frac{\p p}{\p x^i}\frac{\p }{\p \zeta_i}).
\eeq
The integral curves of $H_p$ in $\Sigma_P$ are called null bicharacteristics. Sometimes it is convenient to view these curves on the Lorentzian manifold $(\mbr^4, \tilde g = -dt^2 + g)$. Then the set $\Sigma_P$ consists of light-like vectors of $\tilde g$ and the projections of null bicharacteristics to $\mbr^4$ are light-like geodesics. 

Let $\diag = \{(x, x')\in \mbr^4 \times \mbr^4 : x = x'\}$ be the diagonal of the product manifold and 
\beq
N^*\diag = \{(x, \zeta, x', \zeta')\in T^*(\mbr^4 \times \mbr^4)\backslash 0: x = x', \zeta' = -\zeta\}
\eeq 
be the conormal bundle of $\diag$ minus the zero section. By abuse of notations, we let $\Sigma_P = \{(x, \zeta, x', \zeta')\in T^*\mbr^4\times T^*\mbr^4: p(x, \zeta)= p(x', \zeta') = 0\}$. Then we define $\La_P$ to be the Lagrangian submanifold of $T^*(\mbr^4\times \mbr^4) $ obtained by flowing out $N^*\diag\cap \Sigma_P$ under $H_p$. Here, we also regarded $p(z, \zeta)$ as a function on the product manifold $T^*(\mbr^4\times \mbr^4).$

For two Lagrangian submanifolds $\La_0, \La_1 \subset T^*(\mbr^4\times \mbr^4)$ intersecting cleanly at a co-dimension $k$ submanifold $\Omega \doteq \La_0\cap \La_1$, the space of paired Lagrangian distributions associated with $(\La_0, \La_1)$ is denoted by $I^{p, l}(\La_0, \La_1)$, see \cite{DUV, MU, GrU93} for details. A useful fact is that for $u\in I^{p, l}(\La_0, \La_1)$, we have $u\in I^{p+l}(\La_0\backslash \Omega)$ and $u\in I^p(\La_1\backslash \Omega)$ as Lagrangian distributions which is  recalled in the next paragraph. We know from the results of Melrose-Uhlmann \cite{MU} that the Schwartz kernel of the causal inverse $Q = P^{-1}$ is a paired Lagrangian distribution in $I^{-\frac 32, -\ha}(N^*\diag, \La_P)$. From \cite[Prop.\ 5.6]{DUV}, we also know that $Q: H_{\comp}^{m}(\mbr^4)\rightarrow H^{m+1}_{\loc}(\mbr^4)$ is continuous for $m\in \mbr$.

Let $\La$ be a smooth conic Lagrangian submanifold of $T^*\mbr^4\backslash 0$. Following the standard notation, we denote by $I^\mu(\La)$ the space of Lagrangian distributions of order $\mu$ associated with $\La$, see \cite[Definition 25.1.1]{Ho4}. Such distributions can be represented locally as follows. For $U$ open in $X$, let $\phi(x, \xi): U\times \mbr^N \rightarrow \mbr$ be a smooth non-degenerate phase function that locally parametrizes $\La$ that is, \ $\{(x, d_x\phi): x\in U, d_\xi \phi = 0\} \subset \La.$ Then $u\in I^{\mu}(\La)$ can be locally written as a finite sum of oscillatory integrals
\beq
\int_{\mbr^N} e^{i\phi(x, \xi)} a(x, \xi) d\xi, \ \ a\in S^{\mu + \frac n4 - \frac{N}{2}}(U\times \mbr^N),
\eeq
where $S^\bullet(\bullet)$ denotes the standard symbol class, see \cite[Section 18.1]{Ho3}.  For $u\in I^\mu(\La)$,  the wave front set $\WF(u)\subset \La$ and $u\in H^s(\mbr^4)$ for any $s< -\mu-1$. The principal symbol $\sigma(u)$ of $u \in I^\mu(\La)$ is invariantly defined as a half-density bundle tensored with the Maslov bundle on $\La$, see \cite[Section 25.1]{Ho4}. In local coordinates, these bundles can be trivialized. We remark that we do not specify the order of the principal symbol in the notation but refer to the distribution space for the order.

 A  class of Lagrangian distributions especially important for our purpose is the one of conormal distributions. For a co-dimension $k$ submanifold $Y\subset \mbr^4$, the conormal bundle
\beq
N^*Y = \{(y, \zeta) \in T^*\mbr^4\backslash 0: y\in Y,  \ \  \zeta|_{T_yY} = 0\}
\eeq
is a conic Lagrangian submanifold. The space of conormal distributions to $Y$ of order $\mu$ are denoted by $I^\mu(N^*Y)$. An equivalent definition is that $I^{\mu}(N^*Y)$ consists of $u\in \mcd'(\mbr^4)$ such that 
\beq
L_1 L_2 \cdots L_N u \in {}^\infty H^{\loc}_{-\mu -1}(\mbr^4),
\eeq
where $L_i, i =1, \cdots, N$ are first order differential operators with smooth coefficients tangential to $Y$ and ${}^\infty H^{\loc}_{\bullet}(\mbr^4)$ denotes the Besov space, see \cite[Definition 18.2.6]{Ho3} for details. Such distributions can be represented locally as oscillatory integrals as well.  We know, e.g.\ from \cite[Section 1]{GrU93}, that $I^{\mu}(N^*Y)\subset L^p_{\loc}(\mbr^4)$ for $\mu< -\frac k2 + \frac kp -1$. Examples of conormal distributions are the delta distribution $\delta_Y$ on $Y$, which is in $I^{\frac k2 -1}(N^*Y)$, and a distribution with Heaviside type singularity at $Y$, which is in $I^{-\frac k2 -1}(N^*Y)$.

We restate \cite[Prop.\ 2.1]{GrU93} for the conormal case below.  
\begin{prop}\label{distor}
Let $Y$ be a submanifold of $M$ such that $N^*Y$ intersects $\Sigma_P$ transversally and each bicharacteristics of $P$ intersects $N^*Y$ a finite number of times. For $f\in I^\mu(N^*Y)$,  we have
\beq
v = Q(f)\in I^{\mu-\frac 32, -\ha}(N^*Y, \La_1)
\eeq
where $\La_1 = \La_P\circ N^*Y$ is the flow-out from $N^*Y\cap \Sigma_P$. Furthermore, for $(x, \zeta) \in \La_1\backslash N^*Y$, 
\beq
\sigma(v)(x, \zeta) = \sum_j \sigma(Q)(x, \zeta, y_j, \eta_j) \sigma(f)(y_j, \eta_j)
\eeq 
where $(y_j, \eta_j) \in N^*Y$ is joined to $(x, \zeta)$ by bicharacteristics.  
\end{prop}

\begin{figure}
\centering
\scalebox{0.8}
{
\begin{pspicture}(0,-1.7172384)(9.17,1.7172384)
\definecolor{colour0}{rgb}{0.8,0.8,0.8}
\rput{15.143745}(0.022890234,-0.746676){\psellipse[linecolor=black, linewidth=0.04, fillstyle=solid,fillcolor=colour0, dimen=outer](2.82,-0.28723848)(0.32,0.79)}
\rput{-344.54776}(0.39106882,-1.5369895){\psellipse[linecolor=black, linewidth=0.04, fillstyle=solid,fillcolor=colour0, dimen=inner](5.86,0.67276156)(0.4,1.73)}
\psline[linecolor=black, linewidth=0.04, linestyle=dashed, dash=0.17638889cm 0.10583334cm, dotsize=0.07055555cm 2.0,arrowsize=0.05291667cm 2.0,arrowlength=1.4,arrowinset=0.0]{**->}(0.20776993,-1.101379)(8.93223,1.7269021)
\psline[linecolor=black, linewidth=0.04, linestyle=dashed, dash=0.17638889cm 0.10583334cm, arrowsize=0.05291667cm 2.0,arrowlength=1.4,arrowinset=0.0]{->}(0.26,-1.0772384)(2.52,0.42276153)
\psline[linecolor=black, linewidth=0.04, linestyle=dashed, dash=0.17638889cm 0.10583334cm, arrowsize=0.05291667cm 2.0,arrowlength=1.4,arrowinset=0.0]{->}(0.26,-1.0972384)(2.94,-1.0372385)
\psdots[linecolor=black, dotsize=0.2](2.82,-0.25723848)
\rput[bl](2.68,-1.7172384){$t = t_0$}
\rput[bl](6.08,-1.6572385){$t = t_1$}
\rput[bl](0.0,-1.5572385){$x'$}
\rput[bl](8.9,1.0827615){$\zeta'$}
\end{pspicture}
}
\caption{Distorted plane waves in $\mbr^3$. The two shaded ovals represent the singular support of $f$ at $t = t_0$ and of $v$ at $t = t_1> t_0$.}
\label{figdistor}
\end{figure}
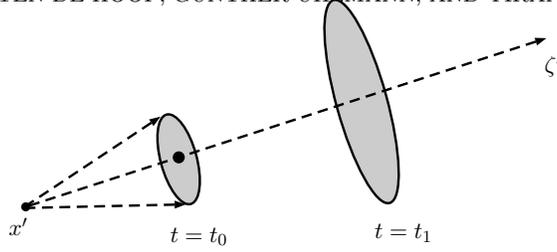

We use the above proposition to construct distorted plane waves. These are generalizations of progressing plane waves but supported near a fixed geodesic. The construction is based on that of \cite{KLU}. For any $(x', \zeta')\in \Sigma_P$, we denote the bicharacteristics from $(x', \zeta')$  by $\Theta_{x', \zeta'}$. Their projections to $\mbr^4$ are denoted by $\gamma_{x', \zeta'}$, which are  light-like geodesics on the Lorentzian manifold $(\mbr^4, \tilde g)$. Here, by abuse of notations, we take $\zeta'$ to be the tangent vector at $x'$ corresponding to $\zeta' \in T_{x'}^*\mbr^4$. This is valid because the non-degenerate metric $g$ induces an isomorphism between $T_{x'}\mbr^4$ and $T^*_{x'}\mbr^4$. For $s_0>0$ a small parameter, we let 
\beq
S(x', \zeta'; s_0) \doteq \{\gamma_{x', \zeta}(\theta)\in \mbr^4: \zeta\in \Sigma_{P, x'}, \|\zeta - \zeta'\| < s_0, \theta > 0\},
\eeq
where the norm is defined using the positive definite metric $\widehat g = dt^2 + g $ on $\mbr^4$.  Notice that as $s_0\rightarrow 0$, $S(x', \zeta'; s_0)$ tends to the geodesic $\gamma_{x', \zeta'}$.  For $t_0 > 0$, we let 
\beqq\label{ymani}
Y(x', \zeta'; t_0, s_0) \doteq S(x', \zeta'; s_0)  \cap \{t = t_0\},
\eeqq
which is a $2$-dimensional surface. See Figure \ref{figdistor}. Then we let 
\beqq\label{lamani}
\La(x', \zeta'; t_0, s_0) \doteq \La_P\circ (N^*S(x', \zeta'; s_0)\cap N^*Y(x', \zeta'; t_0, s_0))
\eeqq
be the flow out. {\em For convenience, we assume that there is no conjugation point on $(\mbr^3, g)$.} We remark that since we essentially consider a local problem in this work, this is not restrictive. Then $S(x', \zeta'; s_0)$ is a co-dimension $1$ submanifold near $\gamma_{x', \zeta'}$ and 
\beq
\La(x', \zeta'; t_0, s_0) = N^*S(x', \zeta'; s_0). 
\eeq
When it is clear from the background, we shall abbreviate the above notations by dropping the dependency on $x', \zeta', t_0, s_0$. 
For $f\in I^\mu(N^*Y)$, using Prop.\ \ref{distor}, we obtain that $v = Qf \in I^{\mu - \frac 32}(\La)$ away from the submanifold $Y$. We conclude that $v$ is conormal to $S$ and we call $v$ a {\em distorted plane wave}.

\section{Local well-posedness of the nonlinear equation}\label{secwell}
For $T> 0$ fixed and $\eps > 0$ small, we consider the well-posedness of the inhomogeneous Cauchy problem
\beq
\begin{gathered}
P u(t, x) + a(t, x) u^2(t, x)  = \eps F(t, x), \ \ (0, T)\times X\\
u(0, x) = \eps f(x), \ \ \p_t u(0, x) = \eps g(x).
\end{gathered}
\eeq
In this section, we use $x\in \mbr^3$ for spatial variables. There is an extensive literature on local and global well-posedness of semilinear wave equations, typically for smooth or power-type nonlinear terms, see e.g.\ Sogge \cite{So}. Here, the problem is that we have a non-smooth nonlinear term. 
If $a(t, x)$ is sufficiently regular, e.g.\ in $H^3(\mbr^4)$ which is an algebra,  it is relatively straightforward to prove the existence for $f, g, F$ sufficiently regular and $\eps$ sufficiently small, see for example \cite[Appendix B]{KLU1}. However, we would like to consider $a(t, x) \in L^\infty(\mbr^4)$ which includes the jump discontinuity. Then the solution is expected to be of only low regularity. We shall give a well-posedness result following  the standard argument using Strichartz type estimates. We remark that we do not intend to pursue the optimal or general result here.

We recall the Strichartz estimates for the Cauchy problem from \cite{MSS} valid for the wave operator on compact Riemannian manifolds without boundary. This is sufficient as we only consider the local problem.  Consider the solution $u$ to the Cauchy problem
\beq
\begin{gathered}
(\p_t^2 + \lap_g)u(t, x) = 0, \ \ (0, T)\times \mbr^3 \\
u(0, x) = f(x), \ \ \p_t u(0, x) = g(x).
\end{gathered}
\eeq
Assume that $f, g$ are supported in a compact set $K\subset \mbr^3$. For $4\leq q < \infty$ and $2\leq r < \infty$, Corollary 3.3 of \cite{MSS} tells that
\beqq\label{eqstr1}
\|u\|_{L^r((0, T); L^q(\mbr^3))} \leq C_T (\|f\|_{H^\gamma(\mbr^3)} + \|g\|_{H^{\gamma-1}(\mbr^3)}), 
\eeqq
with $\gamma = 3(1/2 - 1/q) - 1/r$ and $C_T$ depending on $T > 0$. Here, the norm of the (inhomogeneous) Sobolev spaces are defined by
\beq
\|f\|_{H^\alpha(\mbr^3)} = (2\pi)^{-\frac 32}(\int_{\mbr^3}  (1 + |\xi|^{2})^{\alpha}|\hat f(\xi)|^2 d\xi)^\ha, \ \ \alpha \in \mbr, 
\eeq
where $\hat f$ denotes the Fourier transform of $f$.  Below, we also need the homogeneous Sobolev space $\dot H^\alpha(\mbr^3)$ with norm
\beq
\|f\|_{\dot H^\alpha(\mbr^3)} = (2\pi)^{-\frac 32}(\int_{\mbr^3}   |\xi|^{2\alpha}|\hat f(\xi)|^2 d\xi)^\ha.
\eeq

For our purpose, we shall take $q = r = 4$ in \eqref{eqstr1} so that $\gamma = \ha$. Then we get
\beqq\label{stest1}
\|u\|_{L^4((0, T)\times \mbr^3)} \leq C_T (\|f\|_{H^\ha(\mbr^3)} + \|g\|_{H^{-\ha}(\mbr^3)}). 
\eeqq
It is known that the homogeneous Strichartz estimates imply inhomogeneous estimates from a lemma of Christ and Kiselev \cite{CK}.  Consider
\beq
\begin{gathered}
(\p_t^2 + \lap_g)u(t, x) = F(t, x),\ \ (0, T)\times \mbr^3\\
u(0, x) = f(x), \ \ u_t(0, x) = g(x),\ \ \text{at } t = 0
\end{gathered}
\eeq
where $f, g, F$ are supported in $K$. From \cite[Theorem 3.2 ]{SS} and \eqref{stest1}, we get
\beq
\|u\|_{L^4((0, T)\times \mbr^3)} \leq C_T (\|f\|_{H^\ha(\mbr^3)} + \|g\|_{H^{-\ha}(\mbr^3)} + \|F\|_{L^{\frac 43}((0, T)\times \mbr^3)}),
\eeq
with $C_T$ a generic constant depending on $T$. Together with the conservation of energy for linear wave equations
\beq
\|u(\cdot, T)\|_{\dot H^\ha(\mbr^3)} + \|\p_t u(\cdot, T)\|_{\dot H^{-\ha}(\mbr^3)} = \|f\|_{\dot H^\ha(\mbr^3)} + \|g\|_{\dot H^{-\ha}(\mbr^3)},
\eeq
we obtain 
\beqq\label{stest2}
\begin{gathered}
\|u\|_{L^4((0, T)\times \mbr^3)} + \|u(\cdot, T)\|_{\dot H^\ha(\mbr^3)} + \|\p_t u(\cdot, T)\|_{\dot H^{-\ha}(\mbr^3)}\\
 \leq C_T (\|f\|_{H^\ha(\mbr^3)} + \|g\|_{H^{-\ha}(\mbr^3)} + \|F\|_{L^{\frac 43}((0, T)\times \mbr^3)}).
 \end{gathered}
\eeqq

\begin{prop}\label{reg}
Suppose that $f(x)\in H^{\ha}(\mbr^3), g(x) \in H^{-\ha}(\mbr^3)$, $F(t, x)\in L^{\frac 43}((0, T)\times \mbr^3)$ are supported in $x\in K\subset\subset \mbr^3$. Consider the   Cauchy problem
\beqq\label{eqcauchy}
\begin{gathered}
P u(t, x) + a(t, x) u^2(t, x)  = \eps F(t, x), \ \ (0, T)\times \mbr^3\\
u(0, x) = \eps f(x), \ \ \p_t u(0, x) = \eps g(x),
\end{gathered}
\eeqq
where $a\in L^\infty((0, T)\times \mbr^3), \eps \geq 0$. For $T> 0$ fixed, there exists $\eps_0 > 0$  so that for $\eps \in[0, \eps_0)$, there is a unique solution $u$  such that
\beq
(u, \p_t u) \in C^0((0, T); \dot H^{\ha}(\mbr^3)\times \dot H^{-\ha}(\mbr^3)) \text{ and } u\in L^4((0, T)\times \mbr^3).
\eeq
Moreover, there exists a constant $C$ depending on  $K, T$ such that   
\beq
\|u\|_{L^4((0, T)\times \mbr^3)} \leq C \eps (\|f\|_{H^\ha(\mbr^3)} + \|g\|_{H^{-\ha}(\mbr^3)} + \|F\|_{L^{\frac 43}((0, T)\times \mbr^3)}). 
\eeq
\end{prop}
For later reference, we shall denote the solution space by 
\beqq\label{defx}
\mcx \doteq \{f\in L^4((0, T)\times \mbr^3):  (f, \p_t f) \in C^0((0, T); \dot H^{\ha}(\mbr^3)\times \dot H^{-\ha}(\mbr^3)) \}.
\eeqq
\bpf
We follow a standard argument in the proof of \cite[Theorem 4.1]{So}. 
Consider the existence part. Let $u_{-1} = 0$. We define a sequence $u_m, m = 0, 1, 2, \cdots$ by
\beqq\label{eqite1}
\begin{gathered}
P u_m(t, x)  + a(t, x) u^2_{m-1}(t, x)= \eps F(t, x), \ \ (0, T)\times \mbr^3\\
u_m(0, x) = \eps f(x), \ \ \p_t u_m(0, x) = \eps g(x). 
\end{gathered}
\eeqq
It follows from the finite speed of propagation that all $u_m$ are compactly supported in   $(0, T)\times\mbr^3$. Let
\beq
A_m(T) = \|u_m\|_{L^4((0, T)\times \mbr^3)}, \ \ B_m(T) = \|u_m - u_{m-1}\|_{L^4((0, T)\times \mbr^3)}.
\eeq
We claim that there exists $\eps_0 > 0$ so that 
\beq
A_m(T)\leq 2A_0(T), \ \ B_{m+1}(T) \leq \ha B_m(T) \text{ if } 2A_0(T)\leq \eps_0.
\eeq

For $m, j = 0, 1, 2, \cdots,$ we obtain from \eqref{eqite1} that 
\beqq\label{eqite1}
\begin{gathered}
P (u_m(t, x)-u_{j}(t, x))  + a(t, x) [u^2_{m-1}(t, x) - u^2_{j-1}(t, x)] = 0, \ \ (0, T)\times \mbr^3\\
u_m(0, x)-u_{j}(0, x) = 0, \ \ \p_t [u_m(0, x)-u_{j}(0, x)] = 0. 
\end{gathered}
\eeqq
It follows from the Strichartz estimates \eqref{stest2} and H\"older's inequality that
\beq
\begin{split}
\|u_m - u_{j}\|_{L^4((0, T)\times \mbr^3)} &\leq C \|u^2_{m-1} - u^2_{j-1} \|_{L^{\frac 43}((0, T)\times \mbr^3)}\\
&\leq C  \|u_{m-1}+u_{j-1}\|_{L^2((0, T)\times \mbr^3)} \|u_{m-1} - u_{j-1}\|_{L^4((0, T)\times \mbr^3)}\\
&\leq \ha \|u_{m-1} - u_{j-1}\|_{L^4((0, T)\times \mbr^3)},
\end{split}
\eeq
provided $C[\|u_{m-1}\|_{L^2((0, T)\times \mbr^3)}+ \|u_{j-1}\|_{L^2((0, T)\times \mbr^3)}]  \leq \ha$. Hereafter, $C$ denotes a generic constant. 
Suppose that the first part of the claim is true. Using the fact that $u_m$ are compactly supported, we derive 
\beq
\|u_m\|_{L^2((0, T)\times \mbr^3)}\leq C \|u_m\|_{L^4((0, T)\times \mbr^3)} = C A_m(T) \leq 2C A_0(T).
\eeq
If we take $\eps_0 = 1/(4C)$, we proved that $B_m(T)\leq \ha B_{m-1}(T).$  

Next we prove by induction that $A_m(T)\leq 2A_0(T)$. Suppose this is true for $A_k(T), k\leq m-1$. Taking $j = 0$ in \eqref{eqite1}, we obtain the estimate
\beqq\label{equmest}
\|u_m - u_{0}\|_{L^4((0, T)\times \mbr^3)} \leq \ha \|u_{m-1}\|_{L^{4} ((0, T)\times \mbr^3)} \leq  A_0(T). 
\eeqq
It follows easily that $\|u_m\|_{L^4((0, T)\times \mbr^3)}\leq 2A_0(T)$. This completes the proof of the claim. 

Now we show that the sequence $u_m$ converges to $u$ in  $L^4((0, T)\times \mbr^3)$. From the Strichartz estimates for $u_0$ 
\beqq\label{equmest1}
\|u_0\|_{L^4((0, T)\times \mbr^3)}   \leq C_T \eps (\|f\|_{H^\ha(\mbr^3)} + \|g\|_{H^{-\ha}(\mbr^3)} + \|F\|_{L^{\frac 43}((0, T)\times \mbr^3)}),
\eeqq 
we can choose $\eps = \eps_0$ to satisfy the requirement in the claim. Then it follows that $u_m$ converges to some $u$  in $L^4$, hence in the sense of distribution. Next, it is straightforward to see that 
\beq
\begin{gathered}
\|au^2_{m} - au^2_{m-1}\|_{L^{\frac 43}((0, T)\times \mbr^3)} \leq C\|u_m+ u_{m-1}\|_{L^2((0, T)\times \mbr^3)}  \|u_m- u_{m-1}\|_{L^4((0, T)\times \mbr^3)} \\
\leq C \eps_0\|u_m- u_{m-1}\|_{L^4((0, T)\times \mbr^3)}  \leq C \eps_0 2^{-m}.
\end{gathered}
\eeq
Thus $au^2_m$ converges to $au^2$ in $L^{\frac 43}$ hence also in the sense of distribution. Thus we proved that $u \in L^4((0, T)\times \mbr^3)$ is a weak solution to the Cauchy problem \eqref{eqcauchy}. It follows from  \eqref{equmest} and \eqref{equmest1} that
\beq 
\begin{gathered}
\|u_m\|_{L^4((0, T)\times \mbr^3)}   \leq C_T \eps (\|f\|_{H^\ha(\mbr^3)} + \|g\|_{H^{-\ha}(\mbr^3)} + \|F\|_{L^{\frac 43}((0, T)\times \mbr^3)}) 
 \end{gathered}
\eeq 
for all $m\geq 1$, so the estimates for $\|u\|_{L^4((0, T)\times \mbr^3)}$ follows. 

For the regularity of $u$, observe that for $f, g \in C_0^\infty(\mbr^3)$, the $u_m$ defined in \eqref{eqite1} are all smooth and compactly supported. We can slightly modify the argument for the existence part to show that $(u_m, \p_t u_m)$ is a Cauchy sequence in $C^0((0, T); \dot H^{\ha}(\mbr^3)\times \dot H^{-\ha}(\mbr^3))$ converging to $(u, \p_t u)\in C^0((0, T); \dot H^{\ha}(\mbr^3)\times \dot H^{-\ha}(\mbr^3))$. Finally, for $f \in H^{\ha}(\mbr^3), g \in H^{-\ha}(\mbr^3)$, we use approximation by compactly supported functions to conclude that the solution $u\in \mcx$.

At last, consider the uniqueness of the solution. Suppose that $u, w$ are two solutions and let $U = u-w$. Then we have
\beq
\begin{gathered}
P U  + a(t, x) (u+w) U = 0, \ \ (0, T)\times \mbr^3\\
U = 0, \ \ \p_t U = 0. 
\end{gathered}
\eeq
The Strichartz estimates \eqref{stest2} imply that
\beq
\begin{gathered}
\|U\|_{L^4((0, T)\times \mbr^3)}  \leq C (\|u\|_{L^2}+\|v\|_{L^2})\|U\|_{L^4((0, T)\times \mbr^3)} \leq C\eps_0 \|U\|_{L^4((0, T)\times \mbr^3)}.
 \end{gathered}
\eeq
If $\|U\|_{L^4((0, T)\times \mbr^3)}\neq 0$, we reach a  contradiction when $\eps_0$ is sufficiently small. Thus the solution is unique in $L^4((0, T)\times \mbr^3)$. 
\epf

\section{The nonlinear responses}\label{sec-asymp}
It is convenient to work with a more general setup which includes both the source problem and the Cauchy problem. 
We consider the semilinear wave equation 
\beqq\label{eqsem1}
\begin{gathered}
P(t, x)u  + a(t, x)u^2   =  0, \ \ \text{in } (0, T)\times \mbr^3
\end{gathered}
\eeqq
where $a\in I^{\mu_0}(N^*S_0)\cap L^\infty(\mbr^4)$ for a co-dimension one submanifold $S_0$ of $\mbr^4$ not characteristic for $P$. We assume that $u = u(\eps; t, x) \in C^\infty((0,\eps_0); \mcx)$ is a smooth family of solutions to \eqref{eqsem1} and that $u$ possesses the following asymptotic expansion
\beqq\label{eqasmp2}
u = \eps v + \eps^2 w + o(\eps^2),
\eeqq
where the $o(\eps^2)$ term is in $L^4((0, T)\times \mbr^3)$. We shall call $v$ the linear response and $w$ the nonlinear response. We assume that the linearized solution $v = v_1 + v_2$ where $v_i$ satisfies $Pv_i = 0$ and $v_i \in I^{\mu_i}(N^*S_i), \mu_i < -1,  i = 1, 2$ for co-dimension one submanifolds $S_i$ of $\mbr^4$ characteristic for $P$. Finally, we assume that  $S_i$ intersects $S_j,$ $0\leq i< j \leq 2$ transversally at co-dimension $2$ submanifolds $S_{ij}$, namely $T_p S_i + T_p S_j = T_p\mbr^4, \forall p\in S_i\cap S_j. $ Also, $S_{12}$ and $S_0$ intersect at a co-dimension $3$ submanifold $S_{012}\subset \mbr^4$.  Roughly speaking, we assume that the singular supports of $a, v_1, v_2$ intersect at $S_{012}$ in a transversal way. 

\begin{remark}\label{setupsour}
This setup naturally arises from the source problem 
\beq 
\begin{gathered}
P u(t, x) + a(t, x)u^2(t, x)  =  \eps f(t, x), \ \ \text{in } (0, T)\times \mbr^3,\\
u = 0, \ \ (-\infty, 0)\times \mbr^3,
\end{gathered}
\eeq 
with $\eps$ a small parameter and $f$ constructed in Section \ref{sec-lin}. Then the solution $u$ has the expansion in $\eps$ by Prop.\ \ref{reg}. The linearized solution $v = v_1 + v_2$ where $v_i, i = 1, 2$ are distorted plane waves.
\end{remark}
 
From \eqref{eqsem1} and the linearized equation, we derive that
\beq
\begin{gathered}
P (u - \eps v) + a u^2 = 0 \Longrightarrow
 u = \eps v - Q(au^2).
\end{gathered}
\eeq
Using successive approximation, we get 
\beqq\label{eqapp1}
\begin{split}
u &= \eps v - \eps^2 Q(a v^2) + o(\eps^2)\\
 &= \eps(v_1 + v_2) - \eps^2 [Q(av_1^2) +  Q(av_2^2) + 2 Q(av_1 v_2)]+ o(\eps^2).
\end{split}
\eeqq
Here, the $o(\eps^2)$ term is in $L^4((0, T)\times \mbr^3)$ as a consequence of Prop.\ \ref{reg} and the Strichartz estimates for the linearized (wave) equation. We shall analyze the singularities in the nonlinear response $w$, which is a linear combination of  
\beq
X_{1} = Q(av_1^2), \ \ X_2 = Q(av_2^2), \ \ X_{12} =  Q(av_1 v_2).
\eeq
We use some methods in \cite{LUW} to analyze the singularities of these terms in two subsections.

\subsection{Singularities in $X_1, X_2$}
For these two terms, we claim that the waves can be split into transmitted waves and reflected waves, see Figure \ref{figinter3}. We start with 

\begin{lemma}
Let $S$ be a co-dimension one submanifold of $\mbr^4.$ For $v  \in I^{\mu}(N^*S)$ with $\mu < -1$,  we have $v^2\in I^{2\mu + \frac 32}(N^*S)$. 
\end{lemma}
\bpf
For any $p_0 \in N^*S$, we can choose local coordinates $x = (x^i)_{i = 0}^3$ so that $S = \{x^0 = 0\}$ near $p_0$. Let $\xi = (\xi_i)_{i = 0}^3$ be the dual coordinates on $T^*\mbr^3$. We have $N^*S = \{x^0 = 0, \xi_1 = \xi_2 =\xi_3 = 0\}$. Then we can write $v\in I^{\mu}(N^*S)$ near $p_0$ as an oscillatory integral
\beq
v(x) = \int_{\mbr} e^{i x^0\xi_0} a(x, \xi_0)d\xi_0
\eeq
with $a(x, \xi_0)\in S^{m}(\mbr^4\times \mbr), m = \mu + \ha$. Therefore, we get
\beq
\begin{gathered}
v^2(x) =  \int_{\mbr} \int_\mbr e^{i x^0\xi_0} e^{i x^0 \eta_0} a(x, \xi_0)a(x, \eta_0) d\xi_0 d\eta_0  = \int_\mbr e^{i x^0\zeta_0} b(x, \zeta_0)d\zeta_0, 
\end{gathered}
\eeq
where $\zeta_0 = \eta_0 + \xi_0$ and 
\beq
b(x, \zeta_0) = \int_\mbr a(x, \xi_0)a(x, \zeta_0 - \xi_0)d\xi_0.
\eeq
Let $\eta =  \xi_0/\langle \zeta_0\rangle$. We have
\beq
\p_x^\alpha \p_{\zeta_0}^\beta b(x, \zeta_0) = \langle\zeta_0\rangle^{2m+1-|\beta|} \sum_{\alpha_0 + \alpha_1 = \alpha} \int_{\mbr}\frac{\p_x^{\alpha_0}a(x, \zeta_0\eta)}{\langle\zeta_0\rangle^m} \cdot \frac{\p_x^{\alpha_1}\p_{\zeta_0}^\beta a(x, \zeta_0 - \zeta_0 \eta)}{\langle\zeta_0\rangle^{m-|\beta|}}d\eta.
\eeq
Since $a$ is a symbol of order $m$, we have $|\p_x^\alpha\p_{\xi_0}^\beta a(x, \xi_0)| \leq C \langle\xi_0\rangle^{m-|\beta|}$. Hereafter, $C$ denotes a generic constant. Thus, we estimate
\beq
 |\p_x^\alpha \p_{\zeta_0}^\beta b(x, \zeta_0)| \leq C \langle\zeta_0\rangle^{2m+1-|\beta|} \int_{\mbr} \frac{\langle\zeta_0\eta\rangle^m}{\langle\zeta_0\rangle^m} \cdot \frac{\langle \zeta_0 -\zeta_0\eta\rangle^{m-|\beta|}}{\langle\zeta_0 \rangle^{m-|\beta|}} d\eta. 
\eeq
For $m < -\ha$, the integrand is bounded by $C \langle \eta\rangle^{2m}$ (uniformly for $\zeta_0$) hence the integral is finite. Thus, we showed that $b(x, \zeta_0) \in S^{2m+1}(\mbr^4\times \mbr)$ which implies $v^2 \in I^{2\mu + \frac 32}(N^*S)$ for $\mu<-1$.  
\epf

In our setup, we shall take $\mu_i<-1$ and obtain $v_i^2\in I^{2\mu_i + \frac 32}(N^*S_i), i = 1, 2$ using the lemma. From standard wave front analysis, e.g.\ \cite[Theorem 1.3.6]{Du}, we obtain that $av_i^2$ is a well-defined distribution and 
\beq
\WF(av_i^2)\subset (N^*S_i+ N^*S_0 ) \cup N^*S_i\cup N^*S_0 = N^*S_{0i} \cup N^*S_i\cup N^*S_0.
\eeq  
Here, we used the transversal intersection assumption to get $N^*S_i+ N^*S_0 = N^*S_{0i}.$ More precisely, we can apply \cite[Lemma 1.1]{GrU93} to get 
\beqq\label{eqavi}
a v_i^2 \in  I^{2\mu_i + \frac 32, \mu_0 + 1}(N^*S_{0i}, N^*S_i) + I^{\mu_0, 2\mu_i+\frac 52}(N^*S_{0i}, N^*S_0).
\eeqq
We note that the orders here have different meanings to those in \cite[Lemma 1.1]{GrU93}. 

Now consider $X_i, i = 1, 2$ and recall that $\WF(Q) \subset N^*\diag \cup \La_P$. Away from the intersections $S_{0i}$, we  have
\beq
\WF(X_i) \subset (\La_P\circ  N^*S_{0i}) \cup (\La_P\circ N^*S_0) \cup N^*S_i.
\eeq
Here, we used the fact that $S_i$ are characteristic for $P$ to get $\La_P\circ N^*S_i = N^*S_i$. Observe that this part of $\WF(X_i)$   corresponds to the   transmitted wave. Next, we know that $N^*S_0\cap \Sigma_P = \emptyset$ because $S_0$ is not characteristic for $P$. So it suffices to consider $\La_i \doteq \La_P\circ  N^*S_{0i}, i = 1, 2$ and describe these Lagrangians. 

For some $p\in S_{0i}$, consider the normal vectors $ (1, \alpha) \in N_p^*S_i $ and $ (s, \beta) \in N_p^*S_0$, where $g^*(\alpha, \alpha) =1$ and $s^2 = g^*(\beta, \beta)\neq 1$. Consider their linear combination
\beq
\zeta = a(1, \alpha) + b(s, \beta) = a(1 + bs/a, \alpha + b/a \beta) \in N^*_pS_{0i},  \ \ a, b\in \mbr\backslash 0.
\eeq
Without loss of generality, we can take $a = 1$ and find $b$ so that $\zeta \in \Sigma_P$ from solving a quadratic equation. Now for the Lorentzian metric $\tilde g$, we have 
\beq
\begin{split}
\tilde g^*(\zeta, (s, \beta)) &= -s(1+bs) + g^*(\alpha + b\beta, \beta) \\
& = -s - bs^2 + g^*(\alpha, \beta)   + b g^*(\beta, \beta)  = \tilde g^*((1, \alpha), (s, \beta))
\end{split}
\eeq
Thus the vector $\zeta$ corresponds to the reflected directions after the interaction at $S_0$. Finally, we conclude that $\WF(X_i)\subset \La_i\cup N^*S_i, i = 1, 2,$ with the transmitted waves on $N^*S_i$ and reflected waves on $\La_i$. 

Away from $N^*S_0$ and $N^*S_i$, we obtain from \eqref{eqavi} that $av_i^2 \in I^{\mu_0 + 2\mu_i + \frac 52}(N^*S_{012})$. Therefore, using \cite[Prop.\ 2.1]{GrU93} and wave front analysis, we know that away from $N^*S_0\cup N^*S_i$,
\beq
X_i = Q(av_i^2) \in  I^{\mu_0 + 2\mu_i + 1, -\ha}(N^*S_{0i},  \La_i).
\eeq
Thus $X_i \in I^{\mu_0 + 2\mu_i + 1}(\La_i)$ away from $N^*S_0\cup N^*S_i\cup N^*S_{0i}$ and this is the reflected wave in the nonlinear responses.

\subsection{Singularities in $X_{12}$}
The singularities in $X_{12}$ are analyzed in \cite{KLU} and \cite{LUW} when $S_0$ is also characteristic for $P$. In particular, a conic type singularity is generated. We adapt the analysis to $S_0$ not characteristic for $P$. We start with a wave front analysis to locate the singularities of $X_{12}$.

For $v_i\in I^{\mu_i}(N^*S_i), i = 1, 2,$ we can apply \cite[Lemma 1.1]{GrU93}  to get
\beq
v_1v_2 \in I^{\mu_1, \mu_2+1}(N^*S_{12}, N^*S_1) + I^{\mu_2, \mu_1+1}(N^*S_{12}, N^*S_2).
\eeq
By standard wave front analysis, we know that 
\beq
\begin{gathered}
\WF(a v_1 v_2) \subset N^*S_1\cup N^*S_2\cup N^*S_{12} \cup N^*S_0 \cup N^*S_{01}  
\cup N^*S_{02} \cup N^*S_{012},
\end{gathered}
\eeq
where we used $N^*S_{12} + N^*S_0=N^*S_{012}$ as a consequence of the transversal intersection assumptions. Now consider $\WF(X_{12}).$ We already know that $\La_P\circ N^*S_{0i} = \La_i\cup N^*S_i$. Since $S_i$ are characteristic for $P$, the normal vectors in $N^*S_i$ are light-like vectors for $\tilde g$. As $S_1, S_2$ intersect transversally, it is a fact that the  linear combination of two light-like vectors do not give new light like vectors that is, $N^*S_{12}\cap \Sigma_P = N^*S_1\cup N^*S_2.$ Thus it remains to consider $\La \doteq \La_P \circ N^*S_{012}.$ 

We claim that $S_{012}$ must be a space-like curve, namely the tangent vectors to $S_{012}$ are space-like for $\tilde g$. Consider tangent vectors $(a, \theta), a\in \mbr, \theta\in \mbr^3$ to $S_{012}$ at $p$. If $a = 0$, the vector is space-like. Otherwise, one can rescale the vector so it suffices to consider $(1, \theta), \theta \in \mbr^3$. Observe that light-like vectors $(1, \alpha)\in N_p^*S_1, (1, \beta)\in N_p^*S_2, g^*(\alpha, \alpha) = g^*(\beta, \beta)= 1$  are normal to $S_{012}$. So we get 
\beq
-1 + g(\alpha, \theta) = 0
\eeq
where $\alpha$ becomes the corresponding tangent vector in $T_p\mbr^4.$ Since $g(\alpha, \alpha) = 1$, we conclude that $g(\theta, \theta) > 1$ so that either $(1, \theta)$ is space-like or $\theta = \alpha$. The latter is impossible because the same argument tells $\theta = \beta$ but $\alpha, \beta$ are linearly independent. So we conclude that $(1, \theta)$ is space-like. Notice that $N^*S_{012}\cap \Sigma_P \neq \emptyset$, hence $\La$ is non-empty. Away from the intersections $S_{01}, S_{02}, S_{12}$ and $S_{012}$,  we have
\beq
\WF(X_{12}) \subset N^*S_1\cup N^*S_2 \cup \La_1\cup \La_2 \cup \La. 
\eeq 

We summarize the results above and prove the main result of the paper.
\begin{theorem}\label{main}
Let $S_i, i = 1, 2$ be co-dimension one characteristic submanifolds for $P$ intersecting transversally at $S_{12}$. Let $S_0$ be a co-dimension one submanifold of $\mbr^4$ not characteristic for $P$ and $a\in I^{\mu_0}(N^*S_0)\cap L^\infty(\mbr^4)$. Assume that $S_i, i = 1, 2$ intersects $S_0$ transversally at $S_{0i}$. Suppose that $v_i\in I^{\mu_i}(N^*S_i), \mu_i < -1$ are solutions to $Pv_i = 0$. For $\eps \in (0, \eps_0)$, let $u(\eps; t, x)$ be a one parameter family of solutions in $\mcx$ (defined in \eqref{defx}) to the semilinear wave equation
\beq 
\begin{gathered}
P(t, x) u + a(t, x)u^2  =  0, \ \ \text{in } (0, T) \times \mbr^3
\end{gathered}
\eeq 
 and $u = \eps (v_1+v_2) + \eps^2 w + o(\eps^2)$. Assume that $S_{12}$ intersect $S_0$ transversally at $S_{012}$. Then we have the following conclusions way from the intersection sets $S_{01}, S_{02}, S_{12}$ and $S_{012}$ 
\begin{enumerate}
\item $\WF(w) \subset \La_1\cup \La_2 \cup N^*S_1\cup N^*S_2 \cup \La$.
\item Away from $\La_1\cup \La_2 \cup N^*S_1\cup N^*S_2$,  $w \in I^\mu(\La)$ with $\mu = \sum_{i=0}^2 \mu_i + \ha$. 
\item $\La \cap \WF(w)\neq \emptyset$ if the principal symbols $\sigma(v_i)$ and $\sigma(a)$ are non-vanishing at $S_{012}$. 
\end{enumerate}
\end{theorem}

\bpf
(1). The statement summarizes the results we obtained above. 

(2) and (3). It remains to show $w\in I^\mu(\La)$, in particular, to show that $Q(av_1v_2) \in I^\mu(\La)$ because $X_1, X_2$ terms are smooth near $\La$.

By our assumptions on the intersections of $S_i, i = 0, 1, 2$, for any $p \in S_{012}$, we can find local coordinates $x = (x^i)_{i = 0}^3$ such that $S_i = \{x^i = 0\}$ and $S_{012} = \{x^0 = x^1 = x^2 = 0\}$. We use $\zeta = (\zeta_i)_{i = 0}^3$ as the dual variables to $x$. Then we can express for example $N^*S_0 = \{x^0 = 0, \zeta_1 = \zeta_2 = \zeta_3 = 0\}$ and $N^*S_{012} = \{x^0 = x^1=x^2 = 0, \zeta_3 = 0\}$. In this local coordinates, we can write down the conormal distributions  as
\beq
\begin{gathered}
v_1(x) = \int_{\mbr} e^{i x^1\zeta_1} b_1(x, \zeta_1) d\zeta_1, \ \ v_2(x) = \int_{\mbr} e^{i x^2\zeta_2} b_2(x, \zeta_2) d\zeta_2,\\
a(x) =  \int_{\mbr} e^{i x^0\zeta_0} b_0(x, \zeta_0) d\zeta_0,
\end{gathered}
\eeq
where $b_i \in S^{\mu_i + \ha}(\mbr^4 \times \mbr), i = 0, 1, 2$ are standard symbols. Then we have the multiplication
\beq
a(x)v_1(x)v_2(x) = \int_{\mbr^3} e^{i (x^0\zeta_0 + x^1\zeta_1 + x^2\zeta_2)}  b_0(x, \zeta_0)b_1(x, \zeta_1)b_2(x, \zeta_2) d\zeta_0d\zeta_1 d\zeta_2.
\eeq
We denote $c(x, \tilde \zeta) = b_0(x, \zeta_0)b_1(x, \zeta_1)b_2(x, \zeta_2)$ with $\tilde \zeta = (\zeta_0, \zeta_1, \zeta_2)\in \mbr^3$. 

Now we let $\phi(t), t\geq 0$ be a smooth cut-off function such that $\phi(t) = 1$ for $t\geq 1$ and $\phi(t) = 0$ for $t < \ha$. For $\delta > 0$, we define
\beq
\chi_\delta(\tilde \zeta) = \prod_{i = 0}^2 \phi( \frac{|\zeta_i|}{\delta|\tilde \zeta|}).
\eeq
Then  $\chi_\delta$ is supported on $\{\tilde \zeta \in \mbr^3: \delta |\tilde \zeta| \leq 2|\zeta_i|, i = 0, 1, 2\}$. We conclude that $\chi_\delta c$ is a symbol because
\beq
\begin{gathered}
|\p_x^\alpha \p_{\zeta_0}^{\beta_0} \p_{\zeta_1}^{\beta_1} \p_{\zeta_2}^{\beta_2} (\chi_\delta(\tilde \zeta)c(x, \tilde \zeta))| \leq C_{\chi, \delta} (1+ |\zeta_0|)^{\mu_0 + \ha - \beta_0}(1+ |\zeta_1|)^{\mu_1 + \ha - \beta_1}(1+ |\zeta_2|)^{\mu_2 + \ha - \beta_2}\\
\leq C_{\chi, \delta} (1 + |\tilde \zeta|)^{\mu_0+\mu_1+\mu_2 + \frac 32 - |\beta|}
\end{gathered}
\eeq
where we used $\mu_i < -1, i = 1, 2$ and also $\mu_0 < -1$ because $a$ in particular belongs to $L_{\loc}^p(\mbr^4)$ for all $p>0$.  
We split $av_1v_2$ as 
\beqq\label{eqnonlin}
\begin{gathered}
a(x)v_1(x)v_2(x) = \int_{\mbr^3} e^{i (x^0\zeta_0 + x^1\zeta_1 + x^2\zeta_2)}  \chi_\delta(\tilde\zeta)c(x, \tilde\zeta) d\zeta_0d\zeta_1 d\zeta_2 \\
+ \int_{\mbr^3} e^{i (x^0\zeta_0 + x^1\zeta_1 + x^2\zeta_2)}  (1-\chi_\delta(\tilde\zeta)) c(x, \tilde \zeta) d\zeta_0d\zeta_1 d\zeta_2 \doteq U_1 + U_2.
\end{gathered}
\eeqq
Thus near $S_{012}$ and for any $\delta> 0$, $U_1 \in I^\mu(N^*S_{012})$ with $\mu = \sum_{i = 0}^2 \mu_i + 2$ and $U_2$ is a distribution with $\WF(U_2)$ contained in a  $\delta$ neighborhood of $N^*S_1\cup N^*S_2\cup N^*S_0 \cup N^*S_{12}\cup N^*S_{01}\cup N^*S_{02}$. It is clear from the expression that the symbol of $U_1$ is non-vanishing if $b_i, i = 0, 1, 2$ are non-vanishing. Finally,  $w = Q(av_1v_2) = Q(U_1) + Q(U_2)$. By Prop.\ \ref{distor}, we know that $Q(U_1) \in I^{\mu-\frac 32}(\La)$ away from $N^*S_{012}$ and the symbol is non-vanishing on $\La$. For the other piece, we know that $\WF(Q(U_2))$ is contained in a small neighborhood of $\La_1\cup \La_2 \cup N^*S_1\cup N^*S_2$ and $N^*S_{01}\cup N^*S_{02}\cup N^*S_{12}\cup N^*S_{012}$. This finishes the proof. 
\epf 

From the two subsections, we know that the nonlinear responses consist of reflected waves $X_i \in I^{ 2\mu_i + \mu_0+ 1}(\La_i), i = 1, 2$ and the new wave $X_{12} = I^{ \mu_1 + \mu_2 + \mu_0 + \ha}(\La)$.  These can be distinguished in terms of the order of Lagrangian distributions when $\mu_1 - \mu_2\neq \pm \ha$. 

\section{Linear responses versus nonlinear responses}\label{sec-linres}
For equation \eqref{eqsem1}, we have analyzed the singularities in the asymptotic expansion terms in  \eqref{eqasmp2}. Comparing the wave front sets of the linear response $v$ and the nonlinear response $w$, we find that the differences are the reflected waves on $\La_i, i = 1, 2$ and the conic wave on $\La$. In this section, we demonstrate that if the linear properties of the materials are also different across $S_0$, the linear response may also contain  reflected waves, hence the nonlinear responses on $\La_i$ are potentially indistinguishable. For this reason, it is reasonable to think of the new conic wave at $\La$ as the observable nonlinear effect. 

We continue using the notations  in Section \ref{sec-asymp}. We consider a perturbation problem of \eqref{eqsem1}
\beqq\label{eqsem2}
\begin{gathered}
P u(t, x) + \delta q(t, x) u(t, x) + a(t, x)u^2(t, x)  =  0, \ \ \text{in } (0, T)\times \mbr^3,\\
u(t, x) = \eps (u_1(t, x) + u_2(t, x)), \ \ t < 0,
\end{gathered}
\eeqq
where $\eps, \delta > 0$ are two small parameters. For ease of elaboration, we lower the regularity requirements as follows. We assume that $q, a \in I^{\mu_0}(N^*S_0)$  are compactly supported in $t> 0$ with $\mu_0 < -3$ so that $q, a\in H^{s}(\mbr^4), s = -\mu_0 -1> 2$ which is an algebra. We also assume that the incoming waves $u_i \in I^{\mu_i}(N^*S_i), \mu_i < -3$ and $Pu_i = 0$. Thus $u_i\in H^s(\mbr^4)$ as well. 

We remark that the potential $q$ depending on another small parameter simplifies our argument because it allows us to analyze the singularities in the leading term instead of the full solution. In the linear setting  when the metric $g$ has a conormal singularity across a submanifold so that the coefficient of $\lap_g$ has conormal singularities,  de Hoop, Uhlmann and Vasy studied the transmitted and reflected waves carefully in \cite{DUV}. Also, in the backscattering setting when the potential has a conormal singularity, a similar problem is studied by Greenleaf and Uhlmann \cite{GrU93}. However, both papers require quite complicated analysis to clarify the singularities in the full solution. 
 
Under our regularity assumptions, the local well-posedness of equation \eqref{eqsem2} is essentially known, see e.g.\ \cite[Appendix B]{KLU1}. In particular,   for $\eps$ sufficiently small, there is a unique solution $u \in H_{\loc}^s((0, T)\times \mbr^3)$. 
We also have $u = \eps v + \eps^2 w + o(\eps^2)$ where the $o(\eps^2)$ term is small in $H^s$.  Moreover, since the potential depends on $\delta$,  $v$ and $w$ actually have expansions in $\delta$ as well. Our goal is to analyze the wave front sets of the asymptotic terms of $v, w$.
 
\begin{prop}
Consider equation \eqref{eqsem2} with the above assumptions.
For $\delta, \eps > 0$ sufficiently small, there is a unique solution $u\in H^s_{\loc}((0, T)\times \mbr^3)$ which can be written as 
\beq
u = \eps (u_1 + u_2 + \delta V) + \eps^2 \delta W + O(\eps \delta^2) + O(\eps^3),
\eeq
where the remainder terms are in $H^s_{\loc}((0, T)\times \mbr^3)$. Moreover, away from the sets $S_{0},  S_{12}$, we have
\begin{enumerate}
\item $\WF(u_1 + u_2 + \delta V) \subset  N^*S_1\cup \La_1 \cup N^*S_2\cup \La_2. $ 
\item $\WF(W) \subset N^*S_1\cup \La_1 \cup N^*S_2\cup \La_2 \cup \La.$ 
\end{enumerate}

\end{prop}

\bpf
Since $v$ satisfies the linearized equation, we can write $v = v_1 + v_2$ so that 
\beq 
\begin{gathered}
P v_i(t, x) + \delta q(t, x) v_i(t, x)  =  0, \ \ \text{in } (0, T)\times \mbr^3,\\
v_i(t, x) = u_i(t, x), \ \ t < 0.
\end{gathered}
\eeq 
It suffices to analyze the singularities of $v_1$. Let $\bar v = v_1 - u_1$. We get  
\beq
\begin{gathered}
P \bar v(t, x) + \delta q(t, x) \bar v(t, x)  =  -q(t, x)u_1(t, x), \ \ \text{in } (0, T)\times \mbr^3,\\
\bar v(t, x) = 0, \ \ t < 0.
\end{gathered}
\eeq 
Using the causal inverse $Q = P^{-1}$, we get $\bar v + Q(\delta q\bar v) = -Q( \delta q u_1),$
from which we derive 
\beqq\label{eqseries}
\bar v = \sum_{n = 0}^\infty (-1)^n \delta^n (Q M_q)^n Q (-\delta q u_1),
\eeqq
where $M_q$ denotes the operator of multiplication by $q$. Since $q\in H^{s}(\mbr^4), s> 2$,  we know that $M_q: H^k(\mbr^4)\rightarrow H^k(\mbr^4)$ 
is continuous for $0\leq k\leq s$, see e.g.\ \cite[Section 3.2]{HKM}. We recall that  $u_1\in H^s(\mbr^4)$ and $Q: H^s_{\comp}(\mbr^4)\rightarrow H^{s+1}_{\loc}(\mbr^4)$ is continuous. From the finite speed of propagation (for the linearized equation), we know that each term of \eqref{eqseries} is supported in a compact set of $(0, T)\times \mbr^3$. For $\delta$ sufficiently small, we obtain that the series \eqref{eqseries} converges in $H^{s+1}(\mbr^4)$, and  
\beq
v_1 = u_1 -\delta Q(q u_1) + O(\delta^2),
\eeq
  where the remainder term is in $H^{s+1}(\mbr^4).$

Now we find the singularities in $Q(q u_1)$. Since $q\in I^{\mu_0}(N^*S_0), u_1\in I^{\mu_1}(N^*S_1)$ and $S_0$ intersects $S_1$ transversally, we use \cite[Lemma 1.1]{GrU93} to get 
\beq
qu_1\in I^{\mu_1, \mu_0 + 1}(N^*S_{01}, N^*S_0)  + I^{\mu_0, \mu_1+1}(N^*S_{01}, N^*S_1). 
\eeq
More precisely, we can write $qu_1 = \Phi_1 + \Phi_2$ so that $\Phi_1\in I^{\mu_1, \mu_0 + 1}(N^*S_{01}, N^*S_0)$ microlocally supported away from $N^*S_1$ and $\Phi_2 \in I^{\mu_0, \mu_1+1}(N^*S_{01}, N^*S_1)$ microlocally supported away from $N^*S_0$.  Now consider the action of $Q$ on $qu_1$. Using \cite[Proposition 2.1, 2.2]{GrU93} we obtain that 
\beq
\begin{gathered}
Q(\Phi_2)\in I^{\mu_0 -1, \mu_1}(N^*S_{01}, N^*S_1) +  I^{\mu_0 + \mu_1-\ha, -\ha}(N^*S_{01}, \La_1). 
 \end{gathered}
\eeq
On the other hand, $Q$ acts on $\Phi_1$ as a pseudo-differential operator of order $-2$ so that $Q(\Phi_1) \in I^{\mu_1-2, \mu_0 + 1}(N^*S_{01}, N^*S_0).$
We conclude that the wave front set of $Q(qu_1)$ is contained in $N^*S_{01}\cup N^*S_0\cup N^*S_1\cup \La_1$. The analysis for $v_2$ is the same. So we conclude that 
\beq
v = u_1 + u_2 + \delta V + O(\delta^2)
\eeq
where the wave front set $\WF(V) \subset N^*S_{01}\cup N^*S_0\cup N^*S_1\cup \La_1 \cup N^*S_{02} \cup N^*S_2\cup \La_2. $ Therefore, the linear responses contain reflected and transmitted waves. 

Next, we follows the same lines to analyze the nonlinear response  $w$ which satisfies the equation 
\beq
\begin{gathered}
P w(t, x) + \delta q(t, x) w(t, x)  =  -a(t, x)v^2(t, x), \ \ \text{in } (0, T)\times \mbr^3,\\
w(t, x) = 0, \ \ t < 0.
\end{gathered}
\eeq 
Since $v\in H^{s}(\mbr^4)$ and $a\in H^s(\mbr^4)$, we know that $av^2\in H^s(\mbr^4)$ is well-defined. Similarly, we obtain that 
\beq
w = \sum_{n = 0}^\infty (-1)^n \delta^n (Q M_q)^n Q (-\delta a v^2)
\eeq
which converges in $H^{s+1}(\mbr^4)$ for $\delta$ sufficiently small. So we have 
\beq
w =  \delta W + O(\delta^2), \ \ W = -Q (a (u_1 + u_2)^2).
\eeq 
From wave front analysis as in Section \ref{sec-asymp}, we know that $\WF(W)$ is contained in  
\beq
N^*S_{012} \cup N^*S_{01}\cup N^*S_{02}\cup N^*S_{12} \cup \La_1\cup \La_2 \cup N^*S_1\cup N^*S_2\cup \La.  
\eeq
This completes the proof of the proposition.
\epf

\section{The inverse problem}\label{sec-inv}
As an application of our main results, we address the inverse problem of determining the location of $S_0$ and the principal symbol of $a(t, x)$ using the nonlinear response. We consider a source problem using the construction in Section \ref{sec-lin}.  

\begin{figure}[htbp]
\centering
\scalebox{0.9}
{
\begin{pspicture}(0,-2.3829687)(11.128,2.3829687)
\definecolor{colour0}{rgb}{0.8,0.8,0.8}
\rput{23.652166}(0.15137814,-0.8410245){\psellipse[linecolor=black, linewidth=0.04, fillstyle=solid,fillcolor=colour0, dimen=outer](2.0839999,-0.059031248)(0.256,0.124)}
\rput{-26.645372}(0.27363706,1.0614599){\psellipse[linecolor=black, linewidth=0.04, fillstyle=solid,fillcolor=colour0, dimen=outer](2.378,-0.04703125)(0.25,0.14)}
\psbezier[linecolor=black, linewidth=0.04, fillstyle=solid,fillcolor=colour0](8.493373,1.0001681)(8.173648,0.7095935)(8.201863,0.73798734)(8.198236,0.728333652289474)(8.194611,0.71867996)(8.101056,0.16097312)(8.079479,-0.3613784)(8.057902,-0.88372993)(8.055418,-1.4216963)(8.05585,-1.4223517)(8.056281,-1.423007)(8.012137,-1.3764957)(8.318378,-1.7297602)
\rput{-3.5112748}(0.038724445,0.51533705){\psellipse[linecolor=black, linewidth=0.04, fillstyle=solid,fillcolor=colour0, dimen=outer](8.425824,-0.37402594)(0.11796168,1.376)}
\psdots[linecolor=black, dotsize=0.14](2.208,-0.44703126)
\rput[bl](1.068,2.0729687){$p_1$}
\rput[bl](3.888,2.0529687){$p_2$}
\rput[bl](0.0,-1.2110312){$S_0$}
\rput[bl](2.472,-0.88303125){$p_0$}
\psbezier[linecolor=black, linewidth=0.04](0.568,-1.2870313)(1.4866295,-1.1059985)(1.2244871,-0.7705483)(1.9019738,-0.5071554736024882)(2.5794606,-0.24376266)(3.3950334,-0.86423314)(3.568,-1.2870313)
\rput[bl](3.008,-2.3070312){$\mbr^3$}
\rput{26.82306}(0.80079037,-0.42043075){\psellipse[linecolor=black, linewidth=0.04, fillstyle=solid,fillcolor=colour0, dimen=outer](1.282,1.4689689)(0.202,0.116)}
\rput{-27.360258}(-0.3263845,1.6450335){\psellipse[linecolor=black, linewidth=0.04, fillstyle=solid,fillcolor=colour0, dimen=outer](3.216,1.4929688)(0.228,0.14)}
\rput[bl](0.628,1.3129687){$f_1$}
\rput[bl](3.548,1.2929688){$f_2$}
\psbezier[linecolor=black, linewidth=0.04, linestyle=dashed, dash=0.17638889cm 0.10583334cm, dotsize=0.07055555cm 2.0]{**-}(6.728,2.2729688)(7.2107415,1.7038057)(7.706003,1.001455)(7.9844704,0.3362076163967572)(8.2629385,-0.32903966)(8.592781,-1.0675999)(8.616,-1.8430313)
\psbezier[linecolor=black, linewidth=0.04, linestyle=dashed, dash=0.17638889cm 0.10583334cm, dotsize=0.07055555cm 2.0]{**-}(9.664,2.2089686)(9.177984,1.4804757)(8.992317,1.1809129)(8.657957,0.5952235762548264)(8.323596,0.009534305)(7.9538965,-1.0669491)(7.916,-1.9190313)
\rput[bl](7.088,2.1729689){$p_1$}
\rput[bl](9.868,2.1329687){$p_2$}
\rput[bl](6.216,-1.1310313){$S_0$}
\psbezier[linecolor=black, linewidth=0.04](6.588,-1.2270312)(7.5066295,-1.0459985)(7.292487,-0.6785483)(7.9699736,-0.41515547360248434)(8.647461,-0.15176265)(9.415033,-0.80423313)(9.588,-1.2270312)
\rput[bl](9.368,-2.1030312){$\mbr^3$}
\rput{28.882208}(0.5871043,-0.7418917){\psellipse[linecolor=black, linewidth=0.04, fillstyle=solid,fillcolor=colour0, dimen=outer](1.734,0.76896876)(0.23,0.116)}
\rput{8.615165}(-0.15476413,-0.37140223){\psellipse[linecolor=black, linewidth=0.04, fillstyle=solid,fillcolor=colour0, dimen=outer](2.388,-1.2130313)(0.264,0.13)}
\rput{-28.203058}(-0.02626143,1.4173982){\psellipse[linecolor=black, linewidth=0.04, fillstyle=solid,fillcolor=colour0, dimen=outer](2.8079998,0.7609688)(0.232,0.116)}
\rput{-14.945744}(0.40570244,0.42686898){\psellipse[linecolor=black, linewidth=0.04, fillstyle=solid,fillcolor=colour0, dimen=outer](1.83,-1.3330313)(0.25,0.126)}
\psbezier[linecolor=black, linewidth=0.04, linestyle=dashed, dash=0.17638889cm 0.10583334cm, dotsize=0.07055555cm 2.0,arrowsize=0.05291667cm 2.0,arrowlength=1.4,arrowinset=0.0]{**->}(3.628,2.1329687)(3.1602383,1.4016521)(2.9815447,1.1009282)(2.6597419,0.51296875)(2.3379393,-0.07499064)(1.8358942,-1.1552744)(1.608,-2.0910313)
\psbezier[linecolor=black, linewidth=0.04, linestyle=dashed, dash=0.17638889cm 0.10583334cm, dotsize=0.07055555cm 2.0,arrowsize=0.05291667cm 2.0,arrowlength=1.4,arrowinset=0.0]{**->}(0.708,2.2129688)(1.1907415,1.6438056)(1.686003,0.9414549)(1.9644706,0.2762076163967572)(2.2429383,-0.38903967)(2.4127812,-1.0315999)(2.532,-1.8070313)
\psline[linecolor=black, linewidth=0.04, arrowsize=0.05291667cm 2.0,arrowlength=1.4,arrowinset=0.0]{->}(8.338734,0.9606597)(8.381267,1.3612778)
\psline[linecolor=black, linewidth=0.04, arrowsize=0.05291667cm 2.0,arrowlength=1.4,arrowinset=0.0]{->}(8.174661,-1.7134951)(8.12934,-2.0925674)
\rput[bl](8.832,-0.14303125){$\Lambda$}
\end{pspicture}
}
\label{figa}
\caption{Illustration of Theorem \ref{main1}.  The left picture  shows the setup of the theorem. The two top ovals represent the singular supports of the sources $f_1, f_2$ and the rest represent  the singular supports of distorted plane waves $v_1, v_2$. Similar to Figure \ref{figdistor}, they propagate and concentrate along the geodesics from $p_1, p_2$ (dashed curves). The right picture shows the nonlinear response on $\La$ after the nonlinear interactions at $p_0$, ignoring the transmitted and reflected waves. }
\end{figure}
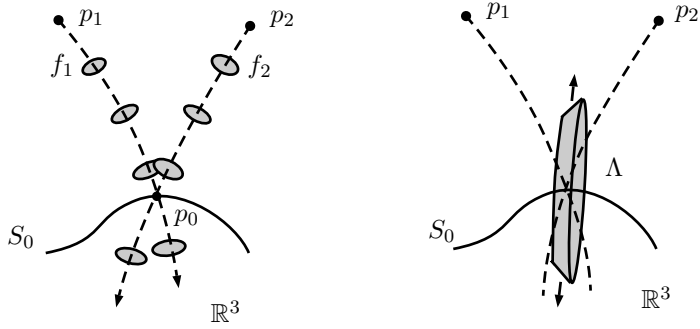

We take two points $(p_i, \zeta_i)\in \Sigma_P, i = 1, 2$  such that the corresponding geodesics $\gamma_{p_i, \zeta_i}$ for the Lorentzian metric $\tilde g = -dt^2 + g$ intersect at $p_0 \in \mbr^4$.   See the left picture of Figure \ref{figa}. For $s_0, t_0 > 0$, let $f_i \in I^{\mu_i+\frac32}(N^*Y_i(p_i, \zeta_i, s_0, t_0))$  and $v_i \in I^{\mu_i}(N^*S_i(p_i, \zeta_i, s_0)), i = 1, 2$ be constructed  as in Section \ref{sec-lin}. Let $S_0$ be a co-dimension one submanifold of $\mbr^4$  not characteristic for $P$, and $a\in I^{\mu_0}(N^*S_0)\cap L^\infty(\mbr^4)$. As in Section \ref{sec-asymp}, we suppose that $S_0, S_1, S_2$ intersect in a transversal way when they intersect. We use the notations $\La_1, \La_2$ in Section \ref{sec-asymp} to denote the Lagrangian submanifolds carrying the reflected waves. Their projections to $\mbr^4$ are denoted by $\widehat S_1, \widehat S_2$ respectively. We denote $\mcs \doteq (\bigcup_{i = 0}^2 S_i) \cup \widehat S_1\cup \widehat S_2$. In particular, we know that this set contains the singular supports of the reflected and transmitted waves in the nonlinear response.

For fixed $T>0$ and  $\eps_1, \eps_2 \in (0, \eps_0)$, we consider the following source problem
\beqq\label{eqsour}
\begin{gathered}
P u(t, x)  + a(t, x) u^2(t, x) = \eps_1 f_1+ \eps_2 f_2, \ \ \text{ in }  (-\infty, T) \times \mbr^3,\\
u(t, x) = 0, \ \ \text{ in } (-\infty, 0)\times \mbr^3.
\end{gathered}
\eeqq
We assume that the exponents $\mu_i, i = 0, 1, 2$ and $\eps_0$ are chosen such that  the well-posedness result Theorem \ref{reg} holds for \eqref{eqsour}. The data set we use for the inverse problem is
\beq
\mcd_{a}(f_1, f_2)\doteq \{u(\eps_1, \eps_2): \text{$u(\eps_1, \eps_2)\in \mcx$ is the unique solution to \eqref{eqsour} for $\eps_1, \eps_2 \in(0, \eps_0)$}\}.
\eeq
We remark that the data set depends on the choice of $(p_i, \zeta_i)$ and $f_i, i = 1, 2$. However, once they are chosen, the data set is a two parameter family of solutions to \eqref{eqsour}.   

\begin{theorem}\label{main1} 
 Suppose that the principal symbols $\sigma(f_i)\neq 0, i = 1, 2$ on $\gamma_{p_i, \zeta_i}$, respectively. Under the above assumptions,  we have
\begin{enumerate}
\item $p_0 \in S_0$ if and only if $\p_{\eps_1}\p_{\eps_2} u(\eps_1, \eps_2)|_{\eps_1 = \eps_2 = 0}$ is not smooth away from $\mcs$ for all $s_0$ small. 
\item If $p_0 \in S_0$, the principal symbol $\sigma(a)$ at $p_0$ is uniquely determined by $\mcd_{a}(f_1, f_2)$. More precisely, suppose $u^{(i)}$ are solutions to \eqref{eqsour} with $a^{(i)} \in I^{\mu_0^{(i)}}(N^*S_0), i = 1, 2$. If $u^{(1)}(\eps_1, \eps_2) = u^{(2)}(\eps_1, \eps_2)$ on $\La$, then the orders $\mu_0^{(1)} = \mu_0^{(2)}$ and the principal symbols $\sigma(a^{(1)}) = \sigma(a^{(2)})$ at $(p_0, \xi_0)\in N^*S_0$. 
\end{enumerate}
\end{theorem}

\bpf
(1). We observed in Remark \ref{setupsour} that the source problem \eqref{eqsour} can be reduced to the setup of Theorem \ref{main}. Following the successive approximation in Section \ref{sec-asymp}, we obtain that 
\beq
\p_{\eps_1}\p_{\eps_2} u(\eps_1, \eps_2)|_{\eps_1 = \eps_2 = 0} = -2Q(a(x)v_1(x)v_2(x)).
\eeq
So the conclusion follows from Theorem \ref{main} when $S_0, S_1, S_2$ intersect at $p_0$. If they do not intersect, the wave front analysis in Section \ref{sec-asymp}  shows that $\WF(Q(a(x)v_1(x)v_2(x)))$ is contained in $(\bigcup_{i = 0}^2 N^*S_i) \cup N^*S_{12}\cup \La_1\cup \La_2$ hence the term is smooth away from the set $\mcs$.

(2). If $\sigma(f_j) \neq 0, j = 1, 2$ on $\gamma_{p_j, \zeta_j}$, we know from Prop.\ \ref{distor} that $\sigma(v_j) \neq 0$ at $(p_0, \xi_j)\in \Sigma_P$. Also, if $u^{(1)}(\eps_1, \eps_2) = u^{(2)}(\eps_1, \eps_2)$  on $\La$, we know from Theorem \ref{main} that $\mcu^{(i)}\doteq \p_{\eps_1}\p_{\eps_2} u^{(i)}(\eps_1, \eps_2)|_{\eps_1 = \eps_2 = 0}, i = 1, 2$ are Lagrangian distributions of the same order on $\La$ away from $\La_1\cup \La_2 \cup N^*S_{12} \cup (\bigcup_{i = 0}^2 N^*S_i)$ with the same principal symbols at $(x, \zeta)\in \La$. By Prop.\ \ref{distor}, we know that  the principal symbols of $\mcu^{(i)}$ at $(p_0, \xi) \in \Sigma_P$ are the same because the matrix $\sigma(Q)(x, \zeta, p_0,  \xi)$ is invertible. In the proof of Theorem \ref{main}, we can read the order and the principal symbols of $\mcu^{(i)}$ at $(p_0, \xi)$ in terms of the principal symbols of $a, v_1, v_2$ at $(p_0, \xi_0), (p_0, \xi_1), (p_0, \xi_2)$ respectively with $ \xi = \sum_{i = 0}^2\xi_i$, see equation \eqref{eqnonlin}. This implies that the order $\mu^{(1)}_0 = \mu^{(2)}_0$ and the principal symbols $\sigma(a^{(1)})(p_0, \xi_0) = \sigma(a^{(2)})(p_0, \xi_0).$
\epf

The nonlinear term can be determined in a special case of piecewise constant functions. The corollary below follows from Theorem \ref{main1} directly. 
\begin{cor}
In addition to the assumptions in Theorem \ref{main1}, we assume that $\Omega$ is a simply connected, bounded open subset of $\mbr^3$ such that $\p \Omega$ is a co-dimension one submanifold of $\mbr^3$. Let $S_0 = \mbr \times \p \Omega$ and $a(t, x) \doteq \alpha \chi_{\Omega}(x), \alpha\in \mbr$, which is conormal  to $S_0$ and in $L^\infty(\mbr^4)$. If $p_0\in S_0$, then  $\alpha$ is uniquely determined by $\mcd_a(f_1, f_2)$. 
\end{cor}


\end{document}